\documentclass[12pt,thmsa]{article}
\usepackage{amsfonts}
\usepackage{amssymb}
\newtheorem{teo}{Theorem}[section]

\newtheorem{lema}[teo]{Lemma}

\newtheorem{prop}[teo]{Proposition}

\newtheorem{obs2}[teo]{Remark}

\newtheorem{tea}{Theorem}[subsection]

\newtheorem{no2}[teo]{Note}

\newtheorem{no3}[tea]{Note}

\newcommand{\Gal}{{\rm Gal}}
\newcommand{\Frob}{{\rm Frob }}
\newcommand{\trace}{{\rm trace}}
\newcommand{\mod}{{\rm mod}}

\newcommand{\Res}{{\rm Res}}

\newcommand{\GL}{{\rm GL}}

\newcommand{\Image}{{\rm Image}}
\newcommand{\PSp}{{\rm PSp}}
\newcommand{\PGSp}{{\rm PGSp}}

\newcommand{\GSp}{{\rm GSp}}
\newcommand{\Sp}{{\rm Sp}}

\newcommand{\cond}{{\rm cond}}
\newcommand{\End}{{\rm End}}
\newcommand{\diag}{{\rm diag}}

\title{ Explicit determination of the images of the
Galois representations attached to abelian surfaces with $\End(A)=\mathbb{Z}$ }

\author{Luis V. Dieulefait  \thanks{supported by TMR - Marie Curie
Fellowship ERBFMBICT983234}\\
Dept. d'·lgebra i Geometria, Universitat de Barcelona\\
Gran Via de les Corts Catalanes 585\\
08007 - Barcelona\\
Spain\\
e-mail: luisd@mat.ub.es}

\begin{document}

\maketitle


\vspace{6mm}





\newpage

\begin{abstract}
We give an effective version of a result reported by
 Serre
asserting that the images of the Galois representations attached to an
abelian surface with $\End(A)= \mathbb{Z}$ are as large as possible for
almost every prime.
Our algorithm depends on the truth of Serre's conjecture for two
dimensional odd irreducible Galois representations. Assuming
this conjecture we determine the finite set of primes with exceptional
image. We also give infinite sets of primes for which we can prove (unconditionally)
that the images of the
corresponding Galois representations are large. We apply the
results to a few examples of abelian surfaces.\\

Mathematics Subject Classification: 11F80, 11G10

\end{abstract}

\newpage

\section{Introduction}

Let $A$ be an abelian surface defined over $\mathbb{Q}$ with
$\End(A) := \End_{\overline{\mathbb{Q}}}(A) = \mathbb{Z}$.\\
Let $ \rho_{\ell}: G_{\mathbb{Q}}\rightarrow \GSp(4, \mathbb{Z}_{\ell} )$ be
the compatible family of Galois representations given by the Galois
action on $ T_{\ell}(A)= A[\ell^{\infty}]( \bar{\mathbb{Q}})$, the Tate
modules of the abelian surface
 (we are assuming that $A$ is
principally polarized). Each $\rho_\ell$ is unramified outside $\ell
N$, where $N$ is the product of the primes of bad reduction of $A$.
If we call $ G_{\ell^{\infty}}$ the image of $ \rho_{\ell}$
then we have the following result of Serre  (cf. [Se 86]):

\begin {teo}
\label {teo:Serre} If $ A$ is an abelian surface over $\mathbb{Q}$ with
$\End(A)=\mathbb{Z}$ and principally polarized, then
$G_{\ell^{\infty}}= \GSp(4, \mathbb{Z}_{\ell} ) $
for almost every $ \ell$.
\end{teo}
Remark: If $ G_{\ell}$ is the image of $\bar{\rho_{\ell}} $, the Galois
 representation
on $\ell$-division points of $ A ( \bar{\mathbb{Q}})$ (and the residual $\mod \; \ell$
representation corresponding to $ \rho_{\ell}$), it is enough to show
that $ G_{\ell} = \GSp(4, \mathbb{F}_{\ell})$ for almost every $\ell$
(cf.  [ Se 86]).\\
Serre proposed the problem of giving an effective version of this
result:\\
``{\it ...partir de courbes de genre $2$ explicites, et  tƒcher de dire … partir
de quand le groupe de Galois correspondant $G_{\ell}$ devient 'gal …
 $ \GSp(4, \mathbb{F}_{\ell}) $ }".
 But Serre's proof depends on certain ineffective results of Faltings
 and therefore does not  solve this problem. \\
 In this article, we will give an algorithm that computes a finite set $
 \mathcal{F}$ of primes containing all those primes (if any)
 with image of the corresponding Galois
 representation exceptional, i.e., different from
 $ \GSp(4, \mathbb{F}_{\ell})$. The validity of our method depends on
 the truth of Serre's conjecture for $2$-dimensional irreducible odd
 Galois representations, conjecture $(3.2.4_{?})$ in [Se 87].
 This means that if $\ell$ is a prime such that
  $ G_{\ell} \neq \GSp(4, \mathbb{F}_{\ell})$ and $\ell \notin
  \mathcal{F}$, then $ G_{\ell}$ has a $2$-dimensional irreducible (odd)
component that violates Serre's conjecture.\\
$ $\\
The method is inspired by the articles [Se 72], [Ri 75], [Ri 85] and [Ri 97]
where   the case of $2$-dimensional Galois
representations is treated.\\
$ $\\
In the examples, we also give infinite sets of primes for which we can
prove the result on the images unconditionally, i.e., without assuming Serre's
conjecture. Results of this kind were previously obtained by Le Duff
only under the extra assumption of semiabelian reduction of the abelian
surface at some prime. Our technique  has two advantages:   it does not
have any restriction on the reduction type of the abelian surface,
 and  in the case of semiabelian reduction it
allows us to prove the
 result on the images (unconditionally) for larger sets of primes.\\
 I want to thank N. Vila, A. Brumer and J-P. Serre for useful remarks
 and comments.

\section{Main Tools}

\subsection{Maximal subgroups of $ \PGSp(4, \mathbb{F}_\ell) $ }

In [Mi 14], Mitchell gives the following classification of maximal
proper subgroups $G$ of $\PSp(4, \mathbb{F}_\ell)$ ($\ell $ odd), as
groups of transformations of the
 projective space having an invariant
linear complex:\\
$ $\\
1) a group having an invariant point and plane\\
2) a group having an invariant parabolic congruence\\
3) a group having an invariant hyperbolic congruence\\
4) a group having an invariant elliptic congruence\\
5) a group having an invariant quadric\\
6) a group having an invariant twisted cubic\\
7) a group $G$ containing a normal elementary abelian subgroup $E$ of
order $16$, with: $ G/E \cong A_5 \; \mbox{or} \;  S_5$ \\
8) a group $G$ isomorphic to $A_6 , S_6$ or $ A_7$.\\
For the relevant definitions see [Hi 85], see also [Bl 17] and [Os 77]
for cases 7) and 8).\\
Remark: This classification is part of a general ``philosophy":
the subgroups of $\GL(n, \mathbb{F}_\ell)$, $\ell$ large, are essentially
 subgroups of Lie type, with some exceptions independent of $\ell$ (see [Se 86]).\\
$ $\\
From this we obtain a classification of maximal proper subgroups $H$ of
$\PGSp(4, \mathbb{F}_\ell)$ with exhaustive determinant. It is similar
to
the above one, except that cases 7) and 8) change according to the
relation between $H$ and $G$, given by the exact sequence:
$$ 1 \rightarrow  G \rightarrow  H \rightarrow  \{  \pm 1 \} \rightarrow 1 .$$

\subsection{The action of Inertia}
From now on we will assume that $\ell$
is a prime of good reduction for the abelian surface
 $A$. Then it follows from results of
Raynaud that the restriction $\bar{\rho}_\ell |_{I_\ell}$ has the
following property (cf. [Ra 74], [Se 72]):
\begin{teo}
\label{teo:Raynaud} If $V$ is a Jordan-H"lder quotient of the
$I_\ell$-module $A[ \ell]( \bar{\mathbb{Q}})$ of dimension $n$ over
$\mathbb{F}_\ell$, then $V$ admits an $\mathbb{F}_{\ell^n}$-vector
space structure of dimension $1$ such that the action of $I_\ell$ on
$V$ is given by a character $ \phi : I_{\ell,t} \rightarrow
\mathbb{F}_{\ell^n}^* $ ($t$ stands for tame) with:
$$ \phi = \phi_1^{d_1}....\phi_n^{d_n}, \quad \qquad (2.1)$$
where $ \phi_i$ are the fundamental characters of level $n$ and $d_i = 0
\; \mbox{or} \; 1$, for every $i=1,2,...,n$.
\end{teo}
All  this statement is proved by Serre in [Se 72] except for the bound for the
 exponents, which
is the  result of Raynaud mentioned above, later generalized by
Fontaine-Messing.\\
$ $\\
We will use the following lemmas repeatedly (cf. [Di 1901]):
\begin{lema}
\label{teo:Dickson} Let $M \in \Sp(4, F)$ be a symplectic transformation
over a field $F$. The roots of the characteristic
polynomial of $M$ can be written as $\alpha , \beta, \alpha^{-1} , \beta^{-1},$
for some $\alpha, \beta$.
\end{lema}
Remark: A similar result holds in general for the groups $\Sp(2n, F)$.\\
In the case of the Galois representations attached to $A$, we know that
$\det(\bar{\rho}_\ell) = \chi^2$, where $\chi$ is the $\mod \; \ell$
cyclotomic character. Therefore, we obtain:
\begin{lema}
\label{teo:Dickson2} The roots of the characteristic
polynomial of $\bar{\rho}_\ell (\Frob \; p) \in G_\ell$ can be written as
 $\alpha , \beta, p/\alpha , p/\beta$
 ($ p \nmid \ell N$).
\end{lema}
Remark: Here $\Frob \;p$ denotes the (arithmetic) Frobenius element, defined up to
 conjugation. The
value of the representation in it is well defined precisely because of
the fact that the representation is unramified at $p$.\\
Proof: Use lemma \ref{teo:Dickson}, $G_\ell \subseteq \GSp(4,
\mathbb{F}_\ell)$, and the exact sequence:
$$ 1 \rightarrow \Sp(4, \mathbb{F}_\ell) \rightarrow \GSp(4, \mathbb{F}_\ell)
\rightarrow \mathbb{F}_\ell^* \rightarrow 1 .$$
$ $\\
Remark: The same is true for $\rho_\ell(\Frob \; p) \in
G_{\ell^\infty}$. Thus, the characteristic polynomial of $\rho_\ell(\Frob \; p)$
has the form:
$$ x^4 - a_p x^3 + b_p x^2 - p a_p x + p^2 $$
with $a_p , b_p \in \mathbb{Z}$, $a_p = \trace (\rho_\ell (\Frob \; p
))$.\\
$ $\\
From formula (2.1), we obtain the following possibilities for
$\bar{\rho}_\ell|_{I_\ell}$ ($ \ell \nmid N$):

$$ \pmatrix{
  1 & * & * & * \cr
  0 & \chi & * & * \cr
  0 & 0 & 1 & * \cr
  0 & 0 & 0 & \chi \cr} ; \quad
\pmatrix{
  \psi_2 & 0 & * & * \cr
  0 &  \psi_2 ^{\ell} & * & * \cr
  0 & 0 & \psi_2 & 0 \cr
  0 & 0 & 0 & \psi_2^\ell \cr} ; $$
 $$ \pmatrix{
   \psi_2 & 0 & * & * \cr
  0 &  \psi_2 ^{\ell} & * & * \cr
  0 & 0 & 1 & * \cr
  0 & 0 & 0 & \chi \cr} ; \quad
  \pmatrix{
  \psi_4^{\ell + \ell^2} & 0 & 0 & 0 \cr
  0 & \psi_4^{\ell^2 + \ell^3} & 0 & 0 \cr
  0 & 0 & \psi_4^{\ell^3 + 1} & 0 \cr
  0 & 0 & 0 & \psi_4^{1 + \ell} \cr} ,$$
$ $\\
where $\psi_i \;$ is a fundamental character of level $i$.\\

\section{Study of the Images}
\subsection{Reducible case: $1$-dimensional constituent}
Let $c$ be the conductor of the compatible family $\{ \rho_\ell \}$.
For the case of the jacobian of a genus $2$ curve, it can be computed
using an algorithm of Q. Liu except for the exponent of $2$ in $c$,
which can easily be bounded using the discriminant of an integral model
of the curve (see [Liu 94]).\\
$ $\\
Suppose that the representation $\bar{\rho}_\ell$ is reducible with a
$1$-dimensional sub(or quotient) representation given by a character
$\mu$. This character is unramified outside $\ell N$ and takes values
in $ \bar{\mathbb{F}}_\ell$, therefore from the description of
$\bar{\rho}_\ell|_{I_\ell}$ given in section 2.2 we have $\mu =
 \varepsilon \chi^i$, with $\varepsilon$ unramified
outside $N$ and $i=0$ or $1$. Clearly $\cond(\varepsilon) \mid c$.
After
semi-simplification, we have:
$$ \bar{\rho}_\ell \cong \varepsilon \chi^i \oplus \pi ,$$
for a $3$-dimensional representation $\pi$ with $\det(\pi) = \varepsilon^{-1}
\chi^{2-i}$.\\
Therefore, $\cond(\varepsilon)^2 \mid c$. Let $d$ be the maximal integer
such that $d^2 \mid c$. If we take a prime $p \equiv 1 \pmod{d}$,
we have $\varepsilon(p)=1$ so  $\chi^i$ is a root of
the characteristic polynomial of $\bar{\rho}_\ell( \Frob \; p)$. This gives:
$$ b_p - a_p (p+1) + p^2 + 1 \equiv 0 \pmod{\ell} , \quad \qquad (3.1)$$
both for $i=0$ and $i=1$ (in agreement with lemma
\ref{teo:Dickson2}).\\
By the Riemann hypothesis, the roots of the characteristic polynomial of
$\rho_\ell ( \Frob \; p )$ have absolute value $ \sqrt{p}$. This gives
automatically bounds for the absolute values of the coefficients
$a_p$ and $b_p$, and from these bounds we see that  for large enough $p$
 congruence (3.1) is not
an equality.
 Therefore, only finitely many primes $\ell$ may verify
(3.1).\\
$ $\\
Variant: Instead of taking a prime $p \equiv 1 \pmod{d}$, we can work
in general with $p$ of  order $f$ in $(\mathbb{Z} / d \mathbb{Z})^*$.
Let $Pol_p(x)$ be the characteristic polynomial of
$\bar{\rho}_\ell (\Frob \; p)$. Then $\varepsilon(p) p^i$ is a root of
$Pol_p(x)$, with $i=0$ or $1$, and $\varepsilon(p)^f = 1 \in
\mathbb{F}_\ell$. Then
$$ \Res(Pol_p(x) , x^f - 1) \equiv 0 \pmod{\ell} \quad \qquad (3.2)$$
where $\Res$ stands for resultant (again, cases $i=0$ and $1$ agree).\\
This variant is used in the examples to avoid computing $Pol_p(x)$ for
large $p$.

\subsection[Reducible case: ``related" $2$-dimensional constituents ]
{Reducible case: ``related" $2$-dimensional \\constituents}

Suppose that, after semi-simplification, $\bar{\rho}_\ell$ decomposes
as the sum of two $2$-dimensional irreducible Galois representations:
$ \bar{\rho}_\ell \cong \pi_1 \oplus \pi_2$. Assume also that these two
constituents are related by lemma \ref{teo:Dickson2}, i.e., if
$\alpha , \beta$ are the roots of the characteristic polynomial of
$\pi_1 (\Frob \; p)$, then $p/\alpha , p/ \beta$ are the roots of that of
$\pi_2 (\Frob \; p)$. If not, then  it follows from
lemma \ref{teo:Dickson2} that  $\alpha =
p/ \beta$, so   $\det(\pi_1) = \det(\pi_2) = \chi$; this case will
 be studied in the next subsection.\\
 Using the description of $\bar{\rho}_\ell|_{I_\ell}$ given in
  section 2.2 we see that one of the following must happen
   (where $\varepsilon$ is a character unramified outside $N$):
 \begin{itemize}
 \item
 Case 1: $\det(\pi_1) = \varepsilon \chi^2 , \quad
  \det(\pi_2) = \varepsilon^{-1} $.
  \item
Case 2: $ \det(\pi_1) = \varepsilon \chi , \quad
   \det(\pi_2) = \varepsilon^{-1} \chi $.
   \end{itemize}
$\bullet$ Case 1: In this case we have the factorization:
$$ Pol_p(x) \equiv (x^2- r x + p^2 \varepsilon(p) ) \;
( x^2 - \frac{r x}{ p \varepsilon(p) } + \varepsilon^{-1}(p)) \quad \pmod{\ell} .$$
As in the previous subsection: $\cond(\varepsilon) \mid d$. Eliminating
$r $ from the equation, we obtain:
$$ Q_p (b_p ,a_p , \varepsilon(p)) := ( \varepsilon(p) b_p -1-p^2 \varepsilon(p)^2 )
(p \varepsilon(p)+1)^2 - a_p^2 p \varepsilon(p)^2 \equiv 0 \pmod{\ell} .$$
If we take $p \equiv 1 \pmod{d}$, we obtain:
$$ (b_p - 1 - p^2) \; ( p+1)^2 \equiv a_p^2 \; p \pmod{\ell} . \quad \qquad (3.3)$$
Again, from the bounds for the coefficients we see that for large enough
$p$ this is not an equality. Thus, only finitely many $\ell$
can satisfy (3.3).\\
$ $\\
Alternatively, for computational purposes, we may take $p$ with
$p^f \equiv 1 \pmod{d}$. Then we have:
$$ \Res( Q_p(b_p , a_p , x) , x^f - 1) \equiv 0 \pmod{\ell} . \quad \qquad (3.4)$$
$\bullet$ Case 2: This case is quite similar to the previous one. We
start with:
$$ Pol_p(x) \equiv (x^2- r x + p \varepsilon(p) ) \;
( x^2 - \frac{r x}{  \varepsilon(p) } + p \varepsilon^{-1}(p)) \quad \pmod{\ell} $$
with $\cond(\varepsilon) \mid d$. From this:
$$ Q'_p (b_p ,a_p , \varepsilon(p)) := ( \varepsilon(p) b_p -p-p \varepsilon(p)^2 )
( \varepsilon(p)+1)^2 - a_p^2  \varepsilon(p)^2 \equiv 0 \pmod{\ell} .$$
Thus, if $p \equiv 1 \pmod{d}$,
$$ 4 ( b_p -2 p ) \equiv a_p^2 \pmod{\ell} . \quad \qquad (3.5) $$
In general, if $p^f \equiv 1 \pmod{d}$,
$$ \Res( Q'_p(b_p , a_p , x) , x^f - 1) \equiv 0 \pmod{\ell} . \quad \qquad (3.6)$$
In this case, the fact that this holds only for finitely many primes $\ell$
is non-trivial. It may be thought of as a consequence of theorem
\ref{teo:Serre}.

\subsection{The remaining reducible case}
As explained above, in the remaining reducible case we have:
 $ \bar{\rho}_\ell^{ss} \cong \pi_1 \oplus \pi_2$ with $\det(\pi_1)=
 \det(\pi_2)=\chi$. We have described in section 2.2 the possibilities for
 $\bar{\rho}_\ell|_{I_\ell}$. This gives for $\pi_1|_{I_\ell}$ and
 $\pi_2|_{I_\ell}$:
$$ \pmatrix{
  1 & *  \cr
  0 &  \chi \cr } \quad   \mbox{or} \quad
\pmatrix{
  \psi_2 & 0 \cr
  0 & \psi_2^{\ell} \cr}.$$
Besides, $\cond(\pi_1) \cond(\pi_2) \mid c $.\\
At this point we invoke Serre's conjecture ($3.2.4_?$) (see [Se 87])
that gives us a control on $\pi_1$ and $\pi_2$. Both representations
should be modular of weight $2$, i.e., there exist two cusp forms
$f_1 , f_2$ with:
$$ \bar{\rho}_{f_1,\ell} \cong \pi_1 , \;  \bar{\rho}_{f_2,\ell} \cong \pi_2 ,
\quad f_1 \in S_2(N_1) , \; f_2 \in S_2(N_2) ,$$
 $N_1 N_2 \mid c$ (we are assuming $\pi_1 , \pi_2$ to be irreducible;
 otherwise,
   they
are covered by section 3.1). Both cusp forms have trivial nebentypus.\\
$ $\\
There are finitely many cusp forms in these finitely many spaces. So
we have an algorithm to detect the primes $\ell$ falling in this case
by comparing characteristic polynomials $\mod \; \ell$, due to
$$ \bar{\rho}_\ell^{ss} \cong \bar{\rho}_{f_1,\ell} \oplus \bar{\rho}_{f_2,\ell} .$$
We take all pairs of integers $N_1 , N_2$ with $N_1 N_2 = c$ and all
pairs of cusp forms $f_1 \in S_2(N_1), \; f_2 \in S_2(N_2)$ (either
newforms or oldforms). If we denote by $Pol_{f_i , p} (x)$ the
characteristic polynomial of $\rho_{f_i , \ell} (\Frob \; p)$
($i=1,2$),
we should have for some such pair $f_1 , f_2$:
$$ Pol_{f_1 , p} (x) Pol_{ f_2 , p} (x) \equiv Pol_p(x) \pmod{\ell}
\quad \qquad (3.7)$$
for every $p \nmid \ell N$.\\
Theorem \ref{teo:Serre} guarantees that this can only happen for
finitely many primes.\\
$ $\\
Remark:
 The Galois representations $\rho_{f_i , \ell}$ attached to
  $f_i$ were constructed by Deligne (cf. [De 71]). The  polynomials
$Pol_{f_i , p} (x)$
are of the form
$$ Pol_{f_i,p}(x)= x^2 - c_p x + p ,$$
where $c_p$ is the eigenvalue of $f_i$ corresponding to the Hecke
operator $T_p$. These eigenvalues, and a fortiori  the characteristic
polynomials $Pol_{f,p}(x)$ for any cusp form $f$, can be computed with
an algorithm of W. Stein (cf. [St]).\\
 The compatible family of Galois
 representations constructed by Deligne,
  in the case of a cusp form $f \in S_2(N)$, shows up in
  the jacobian $J_0(N)$ of the modular curve $X_0(N)$: it agrees
   with a two-dimensional constituent of
   the one attached to the abelian variety $A_f$ corresponding to $f$.

$ $\\
For computational purposes, we introduce the following variant: observe
that either $N_1$ or $N_2$ (say $N_1$) satisfy
 $$ N_1 \mid c , \; N_1 \leq \sqrt{c} . $$
Consider all divisors of $c$ verifying this, maximal (among divisors
of $c$) with this property. Call $S$ the set of such divisors. Then we
are supposing that there exists $f \in S_2(t)$ with $t \in S$ and
$$ \Res(Pol_{f, p}(x) , Pol_p(x)) \equiv 0 \pmod{\ell} ,$$
for every $p \nmid \ell N$. Therefore, for some $t \in S$
$$ \Res (  \prod_{f \in S_2(t)} Pol_{f,p}(x) , Pol_p(x) \;   )
\equiv 0 \pmod{\ell}, \quad \qquad (3.8) $$
for every $p \nmid \ell N$.\\
$ $\\
With this formula we compute in any given example all primes $\ell$
falling in this case.\\
$ $\\
Remark: In all reducible cases (sections 3.1, 3.2 and 3.3) we have
considered reducibility over $\bar{\mathbb{F}}_{\ell}$.

\subsection{Stabilizer of a hyperbolic or elliptic congruence}

If $G_\ell$ corresponds to an
irreducible subgroup inside (its projective image) some of the maximal subgroups in
 cases 3) and 4)  of Mitchell's classification,
 there is a normal subgroup of index $2$ of
 $G_\ell$ such that
$$ 1 \rightarrow M_\ell \rightarrow G_\ell \rightarrow \{ \pm 1 \} \rightarrow 1 ,$$
and the subgroup $M_\ell$ is reducible (not necessarily over $\mathbb{F}_\ell$).\\
In fact, a hyperbolic (elliptic) congruence is composed of all lines
meeting two given skew lines in the projective three dimensional
space over $\mathbb{F}_\ell$ defined over $\mathbb{F}_\ell$
($\mathbb{F}_{\ell^2}$, respectively), called the axes of the congruence
 (see [Hi 85]). The stabilizer of such congruences consists of those transformations that fix
or interchange the two axes, and it contains the normal reducible index two subgroup
 of those transformations that fix both axes.  \\
From the description of $\bar{\rho}_\ell|_{I_\ell}$ given in section 2.2
 we see that if $\ell >
3$, it is contained in $M_\ell$. Therefore, if we take the quotient $G_\ell /M_\ell$
we obtain a representation $ G_\mathbb{Q} \rightarrow C_2$ whose kernel
is a quadratic field unramified outside $N$. Then, there is a quadratic character
$\phi : ( \mathbb{Z}/ c \mathbb{Z})^* \rightarrow C_2$ with:\\
 $\phi(p)=-1 \;  \Rightarrow \; \bar{\rho}_\ell( \Frob \; p)$ is of the
 form:
$$  \pmatrix{
  0 & 0 & * & * \cr
  0 & 0 & * & * \cr
  * & * & 0 & 0 \cr
  * & * & 0 & 0 \cr} .$$
Therefore, $\trace(\bar{\rho}_\ell(\Frob \; p) ) = 0$ , i.e.,
$$ a_p \equiv 0 \pmod{\ell} , \qquad \qquad (3.9) $$
for every $p \nmid \ell N$ with $\phi(p)=-1$.\\
Considering all quadratic characters ramifying only at the primes in
$N$ we detect the primes $\ell$ falling in this case.\\
Once again, from theorem \ref{teo:Serre}, it follows that this set is finite
(of course, this fact strongly depends on the assumption $\End(A)=
\mathbb{Z}$).

\subsection{Stabilizer of a quadric}
This case can be treated exactly as the  one above: assuming again absolute
 irreducibility of the image $G_\ell$, it contains a
normal subgroup of index $2$, and we obtain a quadratic character
unramified outside $N$ verifying formula (3.9). \\
In fact, in this case $\bar{\rho}_\ell$ is the tensor product of two
irreducible $2$-dimensional Galois representations (see [Hi 85], page 28),
 one of them
dihedral (this is the necessary and sufficient condition for the tensor product
 to be symplectic, see [B-R 89], page 51),
 so  the matrices in $G_\ell$ are of the form:
$$  \pmatrix{
  a v & 0 & cv & 0 \cr
  0 & a z & 0 & c z \cr
  b v & 0 & d v & 0 \cr
  0 & b z & 0 & d z \cr} \quad \mbox{or} \quad
    \pmatrix{
  0 & a z & 0 & c z \cr
  a v & 0 & c v & 0 \cr
  0 & b z & 0 & d z \cr
  b v & 0 & d v & 0 \cr} ,$$
depending on the value of the quadratic character $\phi$.

\subsection{Stabilizer of a twisted cubic}
This case is incompatible with the description of $\bar{\rho}_\ell|_{I_\ell} $
 given in section 2.2. In fact, in this case all upper-triangular matrices are of the
 form (see [Hi 85], page 233):
 $$  \pmatrix{
  a^3 & * & * & * \cr
  0 & a^2 d & * & * \cr
  0 & 0 & a d^2 & * \cr
  0 & 0 & 0 & d^3 \cr} .$$
In no case is the subgroup of $G_\ell$ given by  $\bar{\rho}_\ell|_{I_\ell} $
 of this form.

\subsection{Exceptional cases}
The cases already studied cover all possibilities in the classification
except the exceptional groups, i.e., cases 7) and 8).\\
In these cases, comparing the exceptional group
 $H \subseteq \PGSp(4, \mathbb{F}_\ell)$ (its order and structure)
 with the fact that $\mathbb{P}(G_\ell)$
 contains the image of $ \mathbb{P}( \bar{\rho}_\ell|_{I_\ell})$ described in
  section 2.2, we
 end up with the only possibilities ($\ell> 3$):
 $$ \ell = 5 , 7 .$$
 For these two primes, as for any prime we suspect of satisfying $G_\ell \neq
 \GSp(4, \mathbb{F}_\ell) $, we compute several characteristic
 polynomials $Pol_p(x)$ $\mod \; \ell$. At the end, either we prove that
 it must be $G_\ell = \GSp(4, \mathbb{F}_\ell)$ (because the orders of the roots
 of the computed polynomials do not give any other option) or we
 reinforce our suspicion that $\ell$ is exceptional.

 \subsection{Conclusion}
 Having gone through all cases in the classification
  (the stabilizer of a parabolic congruence is reducible, it has an invariant line
   of the complex,
   cf. [Mi 14])
  we conclude that for all primes
  $\ell$ except those whose image, according to our algorithm, may fall
  in a proper subgroup (according to theorem \ref{teo:Serre}, only finitely many)
  the image of $\mathbb{P}(\bar{\rho}_\ell)$ is $\PGSp(4,
  \mathbb{F}_\ell)$.\\
  From this it easily follows that if $\ell$ is not one of
   the finitely many exceptional primes we have
   $G_\ell = \GSp(4, \mathbb{F}_\ell)$
  and applying a lemma of [Se 86] (see also [Se 68]) we obtain
  $G_{\ell^\infty}= \GSp(4, \mathbb{Z}_\ell)$.\\
  Recall that at one step we have assumed the veracity of Serre's
  conjecture ($3.2.4_?$).

  \section{An example}

  We have applied the algorithm to the example given by the jacobian
  of the genus $2$ curve given by the equation:
  $$ y^2 = x^6-x^3-x+1 .$$
  The algorithm of Q. Liu  computes the prime-to-$2$ part of the
  conductor. From this computation and the bound of the conductor in
  terms of the discriminant of an integral equation (cf. [Liu 94]) we
  obtain: $c \mid 2^{12} \cdot 23 \cdot 5$.
 \\
 We exclude a priori
the primes dividing the conductor: $2,5$ and $23$. \\
 We sketch some of the computations performed:\\
 $\bullet$ Reducible cases with $1$-dimensional constituent or two related
 $2$-dimensional constituents:\\
 The maximal possible value of the conductor of $\varepsilon$ is $d=
 64$. We compute the characteristic polynomials of
 $\rho_\ell (\Frob \; p)$ for the primes $p= 229, 257, 641, 769$ and applying
 the algorithm (equations (3.2), (3.4) and (3.6))
 we easily check that no prime $\ell>3$ falls in these
 cases. \\
 Remark: The characteristic polynomials used at this and the remaining
 steps can be found in section 6.\\
$\bullet$ Remaining reducible case:\\
First we describe the set of special divisors of $c$:
$$ S= \{  368, 460, 512 , 640 \} .$$
Then we compute, for each $t \in S$ and each Hecke eigenform $f \in
S_2(t)$,
the characteristic polynomial $Pol_{f,p}(x)$ for
$ p = 3,7,11,13,17,19$  with the
algorithm implemented by W.Stein ([St]).
 Then,
comparing these polynomials with the characteristic polynomials
of $\rho_\ell (\Frob \; p)$ as in formula (3.8) we see that no prime
$\ell > 3$ falls in this case.\\
$\bullet$ Cases ``governed" by a quadratic character:\\
We have to consider all possible quadratic characters $\phi$ unramified
outside
$c$ (there are $15$) and for each of them take a couple of primes
$p$ with $\phi(p)=-1$ and
 $a_p \neq 0$. Applying the algorithm (formula (3.9)) we see
that no prime $\ell > 3$ falls in these cases. At this step we have used the
 values $a_p$ for the primes $p=3,7,13,97,113,569,769$.\\
$\bullet$ Exceptional cases: \\
We compute the reduction of a few characteristic polynomials modulo
$7$ and we find elements whose order (in $\PGSp(4, \mathbb{F}_7)$)
 does not correspond to the
structure of any of the exceptional groups. \\
From all the above computations we conclude:

\begin{teo}
\label{teo:miprimerejemplo} Let $A$ be the jacobian of the genus $2$
curve:
$$ y^2 = x^6 -x^3-x+1 .$$
Let $G_{\ell^\infty}$ be the image of
 $\rho_{A,\ell}  $, the
  Galois representation on $A[ \ell^\infty](\bar{\mathbb{Q}})$, whose conductor
divides $2^{12} \cdot 5 \cdot 23$. Then, assuming Serre's
conjecture ($3.2.4_?$) it holds:
$$ G_{\ell^\infty} = \GSp(4, \mathbb{Z}_\ell)$$
for every $\ell>5 , \; \ell \neq 23$.
\end{teo}
Remark: We are not claiming that the image is not maximal for any of
the four excluded primes.
\section{Unconditional Results and more examples}

\subsection{The case of semiabelian reduction}

  For certain genus $2$ curves one can prove that the image is large
  for an infinite set of primes by using the following results of
  Le Duff [LeD 98]:
  \begin{prop}
  \label{teo:LD1} Let $A$ be an abelian surface defined over $\mathbb{Q}$.
   Suppose that for a prime $p$ of bad reduction of
  $A$, $ \tilde{A}_p^0$ (the connected component of $0$
   in the special fiber of the N'ron Model of $A$ at $p$) is
  an extension of an elliptic curve by a  torus. Then, for every prime
  $\ell \neq p$ with $\ell \nmid \Phi(p)$ (number of connected
  components of $\tilde{A}_p$),
  $G_\ell$ contains a transvection.
  \end{prop}
  Recall that a transvection is an element $u$ such that $\Image(u-1)$ has
   dimension $1$.
  \begin{prop}
  \label{teo:LD2} ([LeD 98]) If $G \subset \Sp(4, \mathbb{F}_\ell)$ is a proper maximal
  subgroup containing a transvection, all its elements have reducible
  (over  $\mathbb{F}_\ell$) characteristic polynomial. Therefore, a transvection
   together with a matrix with  irreducible characteristic
  polynomial
  generate $\Sp(4, \mathbb{F}_\ell)$.
  \end{prop}
  Remark: We can also find in [Mi 14] the list of maximal subgroups of
  $\PSp(4, \mathbb{F}_\ell)$ containing central elations, and a central elation
  is the image in $\PSp(4, \mathbb{F}_\ell)$ of a transvection
   in $\Sp(4, \mathbb{F}_\ell)$. These groups correspond to cases 1) and 3)
   in section 2.1 or to a group having an invariant line of the complex, defined
   over $\mathbb{F}_\ell$.\\
   Recall that $Pol_q(x)$ denotes the characteristic polynomial
    of $\rho_\ell(\Frob \; q)$ for any  prime $q$ of good reduction for
    the abelian surface $A$ and $\ell \neq p$.
  From the two previous results it follows:
  \begin{teo}
  \label{teo:LD3}(Le Duff) Let $p$ be a bad reduction prime
  verifying the condition of proposition \ref{teo:LD1} and $q$ a
  prime with $Pol_q(x)$ irreducible, then for every $\ell \nmid 2 p q \Phi(p)$
  such that $Pol_q(x)$ is irreducible modulo $\ell$,
   $G_\ell = \GSp(4, \mathbb{F}_\ell)$.\\
  If $\Delta_q$ is the discriminant of $Pol_q(x)$ and $\Delta_{Q_q}$ the
  discriminant of $Q_q(x) :=  x^2- a_q x + b_q - 2 q$ the irreducibility
   condition reads:\\
  $ ( \frac{\Delta_q}{\ell} ) = -1 $ and $  ( \frac{\Delta_{Q_q}}{\ell} ) = -1
  $.
  \end{teo}
$ $\\
  EXAMPLE (Le Duff): Take the genus 2 curve:
  $$ C_2 : \quad y^2 = x^5-x+1 . $$
    $A_2=J(C_2)$ has good reduction outside $2,19,151$.
  For $p = 19,151$ the condition in proposition
  \ref{teo:LD1} is satisfied with $\Phi(p)=1$.\\
  Take $q=3$, $Pol_3(x)$ is irreducible and theorem \ref{teo:LD3} gives:
  $G_\ell = \GSp(4,\mathbb{F_\ell})$ for every  $\ell > 3$ with
  $  ( \frac{61}{\ell} ) = -1$ and $  ( \frac{5}{\ell} ) = -1$.\\
  Remark: Of course, considering more irreducible characteristic
  polynomials one can obtain the same result for other  primes.
   In particular, $G_\ell = \GSp(4,\mathbb{F_\ell})$ for $\ell = 19,151$ (cf. [LeD
   98]).\\
  Remark: The example in the previous section also verifies Le Duff's
  condition.\\
$ $\\
  Let us apply our method to this example. The invariants are: \\
  $c = \cond(A_2) \mid 2^8 \cdot 19 \cdot 151$ (computed with
  Liu's algorithm), then $\cond(\varepsilon) \mid d=
  16$; and the set $S= \{  256,604,608 \}$.\\
  In this example, we only have to worry
  about those maximal subgroups in
  Mitchell's classification containing central elations. Therefore, we
   only have to discard the maximal subgroups considered in sections 3.1, 3.3
   and 3.4.\\
   $ $\\
   $\bullet$ The reducible case with 1-dimensional constituent is easily
   handled using the characteristic polynomials (see section 6)
   $ Pol_p(x) $ for $p=17, 97$ and we conclude that no prime $\ell> 2$
    falls in this case.\\
    $\bullet$
Due to the fact that the spaces of modular forms $S_2(t)$ for $t \in S$
 are rather large, we decided to save computations and
  to apply the procedure described in section 3.3, formula (3.8),
   only to the prime
   $p=3$. After computing all resultants of $Pol_3(x)$ with all the
    characteristic polynomials $Pol_{f,3}(x)$ for $f \in S_2(t) , \; t \in
    S$,
    we find the possibly exceptional primes $ \ell > 2 $:
    $$ \ell = 3,5,11,19,29,31,41,61,109,151 .$$
     Having computed
the characteristic polynomials
$$ Pol_p(x) \quad \mbox{for} \quad p=11,41,79,101,199,211 $$
(see section 6) we checked that for each of the ten possibly
exceptional primes $\ell$ listed above one of these six polynomials is irreducible
modulo $\ell$.
Then, applying theorem \ref{teo:LD3}, we conclude that none of
 these primes is exceptional. Thus, no $\ell >2$ has reducible image.\\
 $\bullet$ Cases governed by a quadratic character:
we have to consider all possible quadratic characters $\phi$ unramified
outside
$c$  and for each of them take a couple of primes
$p$ with $\phi(p)=-1$ and
 $a_p \neq 0$. We use the values $a_p$ for
 $p=3,5,97,257$ (see section 6) and an application of the algorithm (formula (3.9))
 proves that the only possibly exceptional primes $\ell>2$ are:
 $$ \ell=3,5,11,97,257 .$$
  We already mentioned that $3,5$ and $11$ are not exceptional.
  Applying  theorem \ref{teo:LD3} again we see that $97$ and $257$ are
  also non-exceptional because $Pol_{11} (x)$ is irreducible modulo
  $97$ and $Pol_{281} (x)$ is irreducible modulo $257$. Then, we
  conclude:
\begin{teo}
\label{teo:ejemplo2a}
    Let $A_2$ be the jacobian of the genus $2$ curve given by the
    equation $y^2 = x^5 - x +1$.
    Assume  Serre's conjecture ($ 3.2.4_? $) (cf. [Se 87]). Then the images
    of the Galois representations on the $\ell$-division points  of $A_2$ are
  $$ G_\ell = \GSp(4,\mathbb{F}_\ell), \quad \mbox{for} \; \mbox{every} \; \ell>2 .$$
  \end{teo}
  Remark: $\bar{\rho}_2$ is also irreducible over $\mathbb{F}_2$. This irreducibility
  for all $\ell$ is equivalent to the fact that $A_2$ is isolated in
  its isogeny class in the sense that any abelian variety
   isogenous to $A_2$ over $\mathbb{Q}$ is isomorphic to $A_2$ over
   $\mathbb{Q}$. Unfortunately, this condition of being isolated is not
   effectively verifiable.\\
$ $\\
  Among the subgroups containing central elations,
  we  have used  Serre's conjecture only to
  eliminate the following one:
  $$ G_\ell \subseteq \{  A \times B \in \GL(2,\mathbb{F}_\ell) \times
  \GL(2,\mathbb{F}_\ell) :
   \det(A)=\det(B)= \chi \}. \qquad (*) $$
  Take $q$ with $Pol_q(x)$ irreducible. If $  ( \frac{\Delta_{Q_q}}{\ell} ) = -1$
   case (*) cannot hold, because the matrices
   $A$ and $B$ would have their traces in $\mathbb{F}_{\ell^2} \smallsetminus
   \mathbb{F}_\ell$. This follows from the factorization:
   $$Pol_q(x) =   \Bigl( x^2 -
    \bigl( \frac { a_q + \sqrt{ \Delta_{Q_q} } } {2} \bigr) x + q \Bigr)
    \Bigl( x^2 -  \bigl( \frac { a_q - \sqrt{ \Delta_{Q_q} } } {2} \bigr) x +
     q \Bigr)  .$$
    Then,  again using $Pol_3(x)$
       we prove
   without using Serre's conjecture  the following
   \begin{teo}
   \label{teo:ejemplo2b}
   The images
    of the Galois representations on the $\ell$-division points of $A_2$ are
  $$ G_\ell =\GSp(4, \mathbb{F}_\ell), \quad \mbox{for every} \quad \ell>3
  \quad  \mbox{with} \quad
    \Bigl( \frac{5}{\ell}  \Bigr) = -1 .$$
    \end{teo}
  Observe that we have obtained an unconditional result that is
  stronger than the one in [LeD 98], because it only uses the condition
  on one of the discriminants (thus, it applies to more primes).
\\
We warn the reader that there is a mistake in [LeD 98], pag. 521,
 the polynomial
$Pol_{11}$ corresponding to this example is wrongly computed. It should read:
$$ x^4+7x^3+31x^2+77x+121.$$

\subsection{Unconditional results in the general case}
 We will show now that even in the case that the condition of
 proposition \ref{teo:LD1} is not verified at any prime,
  we can obtain similar unconditional
 results.\\
 In an arbitrary example, if we do not use Serre's conjecture, there is
 another case to consider (in addition to case (*) ):
 $$ G_\ell \subseteq \{ M \in \GL(2, \mathbb{F}_{\ell^2} ) : \det(M)= \chi \} . \quad (**)$$
 The inclusion of this group in $\GSp(4, \mathbb{F}_\ell)$ is given by
 the map:\\
 $M \rightarrow \diag(M, M^{\Frob})$, where $\Frob$ is the non-trivial
 element in $\Gal( \mathbb{F}_{\ell^2} / \mathbb{F}_\ell)$.\\
 Two tricks allow us to discard this case:\\

 i) Suppose that for a prime $q$, $Pol_q(x)$ decomposes over
 $\mathbb{Q}$ as follows:
 $$ Pol_q(x) = (x^2+Ax+q) (x^2+Bx+q), \quad A \neq B .$$
 Then case (**) cannot hold if $\ell \nmid B-A$ and $\ell \neq q $.\\

 ii) Suppose that $p^{2k+1} \parallel \cond(A)$, then for every
 $\ell \nmid p \Phi(p)$, case (**) cannot hold.
 The condition on $\ell$ is imposed to
 ensure that for these $\ell$  $p^{2k+1} \parallel
 \cond( \bar{\rho}_\ell)$ also holds. \\
$ $\\
 EXAMPLE (Smart): The following  curve is taken from the list
 given in [Sm 97] of all genus $2$ curves defined over $\mathbb{Q}$ with
 good reduction away from $2$:
$$ C_3: \quad y^2 = x(x^4+32x^3+336x^2+1152x-64) ,$$
 $A_3=J(C_3)$, $c \mid 2^{20}$ (this is the
 uniform bound for the $2$-part of the conductor of abelian surfaces over
 $\mathbb{Q}$, cf. [B-K 94]). Le Duff's method cannot be applied to
 this example; the condition of
 proposition \ref{teo:LD1} is not verified at $2$.\\
 We eliminate ALL maximal proper subgroups in Mitchell's classification using
 the characteristic polynomials $Pol_p(x)$ for several
 primes $p$ and $\cond(\varepsilon) \mid 1024$, $S=\{ 1024
 \}$, with the algorithm described in section 3.\\
 To be more precise, the reducible cases treated in section 3.1 and 3.2
 are excluded using the polynomials $Pol_p(x)$ for $p=3,17,19,31$.
Assuming Serre's conjecture, the remaining reducible case is excluded
 using the polynomials $Pol_p(x) $ for $p=7,11,13$.
The cases considered in sections 3.4 and 3.5 are excluded using the
polynomials $Pol_p(x) $ for $p=3,5$. Finally, with the technique described
in section 3.7 we check that $\ell=5,7$ are non-exceptional.\\
All characteristic polynomials used are listed in section 6.\\
 After these computations we  find no  exceptional primes. Then, we conclude
\begin{teo}
\label{teo:ejemplo3a}
    Assume  Serre's conjecture ($ 3.2.4_? $) (cf. [Se 87]). Then the images
    of the Galois representations on the $\ell$-division points of $A_3$ are
  $$ G_\ell = \GSp(4,\mathbb{F}_\ell), \quad \mbox{for} \; \mbox{every} \; \ell>3 .$$
  \end{teo}
$ $\\
 Without Serre's conjecture, trick (i) is used to discard case (**).
 In fact, $Pol_5(x)$ decomposes as in (i) with $A=-2$ and $B= 0$.
  The same happens to $Pol_{17}(x)$.\\
 To deal with case (*) we check that $Pol_3(x)$ is irreducible
 and $\Delta_{Q_3}=12$ (see section 6).
  We obtain the unconditional result:
 \begin{teo}
   \label{teo:ejemplo3b}
   The images
    of the Galois representations on the $\ell$-division points of $A_3$ are
  $$ G_\ell =\GSp(4, \mathbb{F}_\ell), \quad \mbox{for every} \quad \ell>3
  \quad  \mbox{with} \quad
    \Bigl( \frac{3}{\ell}  \Bigr) = -1 .$$
    \end{teo}

  \subsection{Further examples}

  In [Le 91], Lepr'vost gives a genus 2 curve
  over $\mathbb{Q}(t)$ with $13$-rational torsion. For $t=13$ we obtain:
  $$ C : \quad y^2 = -4x^5+300x^4-1404x^3
    +5408x^2-8788x+28561 ,$$
  $A=J(C)$ has $\cond(A)= 2^a \cdot 13^3 \cdot 5^2 \cdot 17^2$.
  Le Duff's condition is not verified at any prime. We can determine
  the image as in the previous example, with or without assuming
  Serre's conjecture (in the
  ``reducible case with 1-dimensional constituent" we  find $\ell = 13$
   as exceptional prime).\\
  Remark: Here trick (ii) eliminates case (**) for every $\ell \neq 13$
  because
  $$ 13^3 \parallel \cond(A) \quad \mbox{and} \quad \Phi(13)= 13 . $$
$ $\\
  Brumer and Kramer (unpublished) have given examples of jacobians of genus $2$
  curves with prime conductor. For them our algorithm determines the
  image with just a few computations. For instance, when applying
  Serre's conjecture  no computation is necessary because
  we have $S= \{ 1 \}$ and $S_2(1) = 0$.\\
  One of these examples is given by the jacobian of the genus two
  curve:
  $$ C: \quad y^2 = x(x^2+1)(1729 x^3 + 45568 x^2 + 25088 x - 76832) .$$
  The conductor of $J(C)$ is $709$.\\
$ $\\
Remark: All the examples of abelian surfaces considered in this article
 verify the
condition $\End(A)= \mathbb{Z}$. This follows in particular from our
result on the images of the attached Galois representations (the condition
on the endomorphism algebra is also necessary for this result to hold).

\section{Computed characteristic polynomials}
We list all the characteristic polynomials $Pol_p(x)$ that have been used in the
examples of the abelian surface $A$ in section 4 and the abelian
surfaces $A_2$ and $A_3$ in section 5.\\
Recall that in any case the polynomial $Pol_p(x)$ is of the form
$$ x^4-a_p x^3 + b_p x^2 - p a_p  x + p^2 $$
so it is enough to give the values $a_p , b_p$.\\

\medskip
\noindent Abelian surface $A$ (section 4):\\

\medskip

\begin{tabular}{|c|c|c|}
  $p$ & $a_p$ & $b_p$ \\
  3 & -3 & 6 \\
  7 & -2 & 6 \\
  11 & -4 & 18 \\
  13 & -5 & 16 \\
  17 & 0 & 22 \\
  19 & -6 & 42 \\
  97 & 6 & 154 \\
  113 & 18 & 250 \\
  229 & 24 & 534 \\
  257 & 15 & 148 \\
  569 & 6 & -118 \\
  641 & 12 & -266 \\
  769 & -6 & 402 \\ \hline
\end{tabular}\\

\newpage

\noindent Abelian surface $A_2 = J(C_2)$ (section 5.1):

\medskip

\begin{tabular}{|c|c|c|}
  $p$ & $a_p$ & $b_p$ \\
  3 & -3 & 7 \\
  5 & -5 & 15 \\
  17 & -3 & 16 \\
  97 & -8 & 86 \\
  257 & -11 & -113 \\
  11 & -7 & 31 \\
  41 & -7 & 72 \\
  79 & 7 & 75 \\
  101 & -8 & -16 \\
  199 & 25 & 338 \\
  211 & -17 & 103 \\
  281 & 1 & 148 \\ \hline
\end{tabular}

\medskip

\medskip
\noindent Abelian surface $A_3 = J(C_3)$ (section 5.2):

\medskip

\begin{tabular}{|c|c|c|}
  $p $ & $a_p$ & $b_p$ \\
  3 & 2 & 4 \\
  5 & 2 & 10 \\
  7 & -2 & 2 \\
  11 & -2 & 12 \\
  13 & -6 & 18 \\
  17 & 4 & 22 \\
  19 & -2 & -4 \\
  31 & 4 & 46 \\ \hline
\end{tabular}

\section{Bibliography}\

[B-R 89] Blasius, D., Ramakrishnan, D., {\it Maass Forms and Galois
Representations}, in ``Galois groups over $\mathbb{Q}$", Y. Ihara, K.
Ribet, J-P. Serre, eds., MSRI Publications, Springer-Verlag (1989)
33-78

[B-K 94]- Brumer, A., Kramer, K., { \it The conductor of an abelian
variety}, Compositio Math. {\bf 92} (1994) 227-248

[B-K]-  \_\_\_\_\_\_\_\_\_, { \it Non-existence of certain semistable abelian
varieties}, preprint

[Bl 17]- Blichfeldt, H., {\it Finite collineation groups}, (1917)
University of Chicago Press

[De 71]- Deligne, P., {\it Formes modulaires et repr'sentations
$\ell$-adiques}, Lect. Notes in Mathematics {\bf 179} (1971) 139-172

[Di 1901]- Dickson, L., {\it Canonical forms of quaternary abelian
substitutions in an arbitrary Galois field}, Trans. Amer. Math. Soc. {\bf 2}
(1901) 103-138

%
[Hi 85]- Hirschfeld, J., {\it Finite projective spaces of three
dimensions}, (1985) Clarendon Press - Oxford

[Le 91]- Lepr'vost, F., {\it Famille de courbes de genre $2$ munies d'une classe
de diviseurs rationnels d'ordre $13$}, C. R. Acad. Sci. Paris { \bf 313
}
S'rie I (1991) 451-454

[LeD 98]- Le Duff, P., {\it Repr'sentations Galoisiennes associ'es aux
points d'ordre $\ell$ des jacobiennes de certaines courbes de genre $2$},
 Bull. Soc. Math. France {\bf 126} (1998) 507-524

[Liu 94]- Liu, Q., {\it Conducteur et discriminant minimal de courbes de genre
$2$}, Compositio Math. {\bf 94} (1994) 51-79

[Mi 14]- Mitchell, H., {\it The subgroups of the quaternary abelian
linear group}, Trans. Amer. Math. Soc. {\bf 15} (1914) 379-396

[Os 77]- Ostrom, T., {\it Collineation groups whose order is prime to the
characteristic}, Math. Z. {\bf 156} (1977) 59-71

[Ra 74]- Raynaud, M., {\it Sch'mas en groupes de type $(p,...,p)$},
Bull. Soc. Math. France {\bf 102} (1974) 241-280

[Ri 75]- Ribet, K.A., {\it On $\ell$-adic representations attached
 to modular forms},
Invent. Math. {\bf 28} (1975)  245-275

[Ri 85]- \_\_\_\_\_\_\_\_\_, {\it On $\ell$-adic representations attached
 to modular
forms II}, Glasgow Math. J. {\bf 27} (1985) 185-194

[Ri 97]- \_\_\_\_\_\_ , {\it Images of semistable Galois
representations}, Pacific J. of Math. {\bf 181} (1997)

[Se 68]- Serre, J-P., {\it Abelian $\ell$-adic representations and elliptic
curves}, (1968) Benjamin

[Se 72]- \_\_\_\_\_\_ ,  {\it Propri't's galoisiennes des points d'ordre
fini des courbes elliptiques}, Invent. Math. {\bf 15} (1972) 259-331

[Se 86]- \_\_\_\_\_\_ , {\it Oeuvres}, vol. 4, 1-55, (2000) Springer-Verlag

[Se 87]- \_\_\_\_\_\_ , {\it Sur les repr'sentations modulaires de degr'
$2$ de $\Gal(\bar{\mathbb{Q}} / \mathbb{Q})$}, Duke Math. J. {\bf 54}
(1987) 179-230

[Sm 97]- Smart, N.P., { \it $S$-unit equations, binary forms and curves of
genus $2$}, Proc. Lond. Math. Soc., III. Ser. 75, No.2, (1997) 271-307

[St]- Stein, W., {Hecke: The Modular Forms Calculator},
   available at:
 http://shimura.math.berkeley.edu/$\sim$was/Tables/hecke.html

\end{document}